\documentclass{article}
\usepackage{amsmath} 
\usepackage{amsfonts}
\usepackage{amsthm}
\usepackage{amssymb}
\usepackage{amscd}
\newtheorem{Thm}{Theorem}
\newtheorem{Def}{Definition}
\newtheorem{Lem}{Lemma}
\newtheorem{Rem}{Remark}
\newtheorem{Prop}{Proposition}

\title{No $C^1$-recurrence of iterations of symplectomorphisms}
\author{Yoshihiro Sugimoto}
\date{}

\begin{document}

\maketitle

\begin{abstract}
    In this article, we study the  behaviour  of iterations of symplectomorphisms and Hamiltonian diffeomorphisms on symplectic manifolds. We prove that symplectomorphisms and Hamiltonian diffeomorphisms do not have $C^1$-recurrence on negatively monotone symplectic manifolds. This is a generalization of the results of the study by Polterovich, Ono, Atallah-Shelukhin. Hamiltonian group actions play very important roles in symplectic topology. We see that negatively monotone symplectic manifolds are far from being Hamiltonina G-maniofolds.
\end{abstract}

\section{Introduction}
In \cite{P}, Polterovich introduced ``growth sequence" of a diffeomorphism $f$ of a smooth compact manifold $M$ as follows:
\begin{equation*}
    \Gamma_n(f)=\max\{\ \max_{x\in M}|d_xf^n|, \ \max_{x\in M}|d_x^{-n}| \ \}, \ n\in \mathbb{N}
\end{equation*}
Here, $|d_xf^n|$ is the operator norm of the differential map
\begin{equation*}
    d_xf^n:T_xM\longrightarrow T_{f^n(x)}M
\end{equation*}
caluculated with respect to some Riemannian metric on $M$. Polterovich studied the behabior of ${\Gamma_n(f)}$ as $n$ goes to $+\infty$ for symplectomorphism $f$. Let ${(M,\omega)}$ be a closed symplectic manifold. We denote the group of symplectomorphisms by ${\textrm{Symp}(M,\omega)}$.
\begin{equation*}
    \textrm{Symp}(M,\omega)=\{\phi\in\textrm{Diff}(M)\ | \ \phi^*\omega=\omega\}
\end{equation*}
$\textrm{Symp}_0(M,\omega)$ is the identity component of ${\textrm{Symp}(M,\omega)}$. In other words, any element of ${\textrm{Symp}_0(M,\omega)}$ can be connected to the identity via an isotopy of symplectomorphisms. 
\begin{gather*}
    \mathrm{Symp}_0(M,\omega)=\bigg\{\phi\in \mathrm{Symp}(M,\omega) \ \bigg| \ \begin{matrix}
        \exists \ \mathrm{smooth \ isotopy} \{\phi^t\}\subset \mathrm{Symp}(M,\omega) \\ \mathrm{s.t.} \ \phi^0=\mathrm{Id}, \phi^1=\phi\end{matrix}\bigg\}
\end{gather*}

Polterovich proved that if ${\pi_2(M)=0}$ and ${\phi\in \textrm{Symp}_0(M,\omega)\backslash \{\textrm{Id}\}}$ has a fixed point of contractible type, 
\begin{equation*}
    \Gamma_n(f)\rightarrow +\infty \ \ (n\rightarrow +\infty)
\end{equation*}
holds. In particular, any Hamiltonian diffeomorphism does not have $C^1$-recurrence property if $\pi_2(M)=0$ holds (\cite{P}). A symplectic manifold ${(M,\omega)}$ is called negatively monotone if there is a negative constant ${\kappa<0}$ such that 
\begin{equation*}
    c_1|_{\pi_2(M)}=\kappa \cdot \omega|_{\pi_2(M)}
\end{equation*}
holds. Here, ${c_1\in H^2(M)}$ is the first Chern class of the tangent bundle of $M$. Ono proved that there is no Hamiltonian $S^1$-action on any negatively monotone symplectic manifold (\cite{O}). Recently, Atallah-Shelukhin proved that there is no Hamiltonian torsion on negatively monotone symplectic manifolds (\cite{AS}). In other words, there is no Hamiltonian diffeomorphism $\phi$ (${\phi\neq \mathrm{Id}}$) such that ${\phi^k=\mathrm{Id}}$ holds for some positive integer ${k>1}$. In this paper, we generalize these results and prove that there is no $C^1$-recurrence of iterations of symplectomorphisms and Hamiltonian diffeomorphisms on negatively monotone symplectic manifolds. In particular, the set ${\{\phi^k\}_{k\in \mathbb{Z}}}$ is discrete in $C^1$-topology. Moreover, ${(M,\omega)}$ is far from being a Hamiltonian $G$-manifold.
 
 For any smooth function ${H\in C^{\infty}(M)}$, we define the Hamiltonian vector field ${X_H}$ by the following relation:
 \begin{gather*}
     \omega(X_H,\cdot)=-dH
 \end{gather*}
We also consider an $S^1$-dependent ($=$ $1$-periodic) Hamiltonian function ${H\in C^{\infty}(S^1\times M)}$ and an $S^1$-dependent Hamiltonian vector field $X_H$ by the same formula. The time $1$ flow of the vector field $X_H$ is called Hamiltonian diffeomorphism generated by $H$. We denote this Hamiltonian diffeomorphism by $\phi_H$. The set of all Hamiltonian diffeomorphisms is called the Hamiltonian diffeomorphism group and we denote it by ${\mathrm{Ham}(M,\omega)}$, i.e.,
\begin{gather*}
    \mathrm{Ham}(M,\omega)=\{\phi_H \ | \ H\in C^{\infty}(S^1\times M)\}
\end{gather*}
${\mathrm{Ham}(M,\omega)}$ is a subgroup of ${\mathrm{Symp}_0(M,\omega)}$. The composition of ${\phi_H, \phi_K}$ is generated by a Hamiltonian function
\begin{gather*}
    H\sharp K(t,x)=H(t,x)+K(t,(\phi_H^t)^{-1}(x)).
\end{gather*}
Moreover, ${\phi_{H\sharp K}^t=\phi_H^t\circ \phi_K^t}$ holds. The inverse of ${\phi_H}$ is generated by a Hamiltonian function
\begin{gather*}
    \overline{H}(t,x)=-H(t,\phi_H^t(x)).
\end{gather*}
This $\overline{H}$ satisfies ${\phi_{\overline{H}}^t=(\phi_H^t)^{-1}}$.
We also consider $k$-th power of Hamiltonian diffeomorphisms for any integer ${k\in \mathbb{N}}$. We define ${H^{(k)}\in C^{\infty}(S^1\times M)}$ as follows:
\begin{gather*}
    H^{(k)}(t,x)=kH(kt,x)
\end{gather*}
It is straightforward to see that ${\phi_{H^{(k)}}=(\phi_H)^k}$ holds. In other words, ${H^{(k)}}$ generates the $k$-th power of $\phi_H$.

The main results of this paper are as  follows:
\begin{Thm}[no $C^1$-recurrence]
Let $(M,\omega)$ be a closed negatively monotone symplectic manifold.
\begin{enumerate}
    \item There is a $C^1$-small open neighborhood ${\mathcal{U}\subset \textrm{Ham}(M,\omega)}$ of the identity such that for any ${\phi \in \textrm{Ham}(M,\omega)\backslash \{\textrm{Id}\}}$, we can choose a positive integer ${N_{\phi}>0}$ so that 
    \begin{equation*}
        k\ge N_{\phi}\Longrightarrow \phi^k\notin \mathcal{U}
    \end{equation*}
    holds.
    \item Assume that the Euler number of $M$ is not zero (${\chi(M)\neq 0}$) or ${\pi_1(M)}$ has finite center. Then, there is a $C^1$-small open neighborhood ${\mathcal{V}\subset \textrm{Symp}_0(M,\omega)}$ of the identity such that for any ${\psi \in \textrm{Symp}_0(M,\omega)\backslash \{\textrm{Id}\}}$, we can choose a positive integer ${M_{\psi}>0}$ so that 
    \begin{equation*}
        k\ge M_{\psi}\Longrightarrow \psi^k\notin \mathcal{V}
    \end{equation*}
    holds.
\end{enumerate}
\end{Thm}

\section{Flux homomorphism and flux group}

 Let ${\widetilde{\textrm{Symp}}_0(M,\omega)}$ be the universal cover of ${\textrm{Symp}_0(M,\omega)}$. A point in ${\widetilde{\textrm{Symp}}_0(M,\omega)}$ is a homotopy class of smooth paths ${\{\phi^t\}\subset \textrm{Symp}(M,\omega)}$ with fixed endpoints ${\phi^0=\textrm{Id}}$ and ${\phi^1=\phi}$. Let ${\{\phi^t\}_{t\in [0,1]}\subset \textrm{Symp}(M,\omega)}$ be a smooth symplectic isotopy generated by a vector field ${X_t}$.
 \begin{equation*}
     \frac{d}{dt}\phi^t(x)=X_t(\phi^t(x))
 \end{equation*}
Note that ${\iota_{X_t}\omega\in \Omega^1(M)}$ is a closed $1$-form because ${\{\phi^t\}}$ is a symplectic isotopy and 
\begin{gather*}
    0=L_{X_t}\omega=d(\iota_{X_t}\omega)+\iota_{X_t}(d\omega)=d(\iota_{X_t}\omega)
    \end{gather*}
 holds. The flux homomorphism of this isotopy is defined as follows:
\begin{equation*}
    \textrm{Flux}(\{\phi^t\})=\Big[\int_0^1\iota_{X_t}\omega dt\Big]\subset H^1(M:\mathbb{R})
\end{equation*}
The flux homomorphism is invariant under any homotopy with fixed endpoints. So, it is defined on the universal cover.
\begin{equation*}
    \textrm{Flux}: \widetilde{\textrm{Symp}_0}(M,\omega)\longrightarrow H^1(M:\mathbb{R})
\end{equation*}
The flux group ${\Gamma\subset H^1(M:\mathbb{R})}$ of a symplectic manifold ${(M,\omega)}$ is the image of the fundamental group of ${\textrm{Symp}_0(M,\omega)}$ under the flux homomorphism.
\begin{equation*}
    \Gamma=\textrm{Flux}(\pi_1(\textrm{Symp}_0(M,\omega)))
\end{equation*}
The most important property of $\Gamma$ is that it is a discrete subgroup of ${H^1(M:\mathbb{R})}$ for any closed symplectic manifold ${(M,\omega)}$(\cite{O}).

\section{Floer homology and mean index}
Let ${(M,\omega)}$ be a $2n$-dimensional closed symplectic manifold and we fix a $1$-periodic Hamiltonian function 
\begin{gather*}
    H:S^1\times M\longrightarrow \mathbb{R}.
\end{gather*}
\begin{Def}
A Hamiltonian function ${H\in C^{\infty}(S^1\times M)}$ is called non-degenerate if the differential map
    \begin{gather*}
        (d\phi_H)_x:T_xM\longrightarrow T_xM
    \end{gather*}
    does not have $1$ as an eigenvalue for any fixed point $x\in \mathrm{Fix}(\phi_H)$. In other words, the graph of $\phi_H$ intersects the diagonal ${\Delta_M\subset M\times M}$ transverally.
\end{Def}
We denote the set of contractible $1$-periodic orbits of ${\phi_H^t}$ by ${\mathcal{P}(H)}$.
\begin{gather*}
    \mathcal{P}(H)=\{x:S^1\rightarrow M \ | \ \dot{x}(t)=X_{H_t}(x(t)), \ x:\mathrm{contractible}\}
\end{gather*}
We define a covering space of ${\mathcal{P}(H)}$ as follows:
\begin{gather*}
    \widetilde{\mathcal{P}}(H)=\{(u,x) \ | \ u:D^2\rightarrow M, x\in \mathcal{P}(H), \partial u=x\}/\sim
\end{gather*}
The equivalence relation $\sim$ is defined as follows:
\begin{gather*}
    (u,x)\sim (v,y)\Longleftrightarrow \begin{cases}  x=y  \\
    \int_{S^2}(u\sharp \overline{v})^*\omega=0 \\  \int_{S^2}(u\sharp \overline{v})^*c_1=0    \end{cases}
\end{gather*}

This ${u\sharp \overline{v}}$ is a sphere map ${u\sharp \overline{v}:S^2\rightarrow M}$ obtained by gluing ${u:D^2\rightarrow M}$ and ${\overline{v}:D^2\rightarrow M}$ along the boundary. 
Here, $\overline{v}$ is the disc with the opposite orientation. ${c_1\in H^2(M)}$ is the first Chern class of ${(M,\omega)}$. We denote the equivariance class of ${(u,x)}$ by ${[u,x]}$. For any ${[u,x]}$, the action functional ${A_H([u,x])}$ is defined by
\begin{gather*}
    A_H([u,x])=-\int_{D^2}u^*\omega+\int_0^1H(t,x(t))dt .
\end{gather*}
We have the Conley-Zehnder index ${\mu_{CZ}}$ of ${\widetilde{\mathcal{P}}(H)}$ (see \cite{SZ,RS}). We normalize $\mu_{CZ}$ so that the Conley-Zehnder index of a local maximum of a $C^2$-small Morse function is equal to $n$. ${\pi_2(M)}$ acts ${\widetilde{\mathcal{P}}(H)}$ naturally and it changes the Conley-Zehnder index as follows:
\begin{gather*}
    \mu_{CZ}([u\sharp A,x])=\mu_{CZ}([u,x])-2c_1(A) \ \ \ (\forall A\in \pi_2(M))
\end{gather*}
We also have the mean index ${\mu([u,x])\in \mathbb{R}}$ (\cite{SZ}). The mean index has the following properties.
\begin{enumerate}
    \item $|\mu([u,x])–\mu_{CZ}([u,x])|\le n$
    \item ${\mu([u^{(k)},x^{(k)}])=k\mu([u,x])}$
    \end{enumerate}
Here ${[u^{(k)},x^{(k)}]\in \widetilde{\mathcal{P}}(H^{(k)})}$ is the natural $k$-th iteration of ${[u,x]}$. The Floer chain complex ${CF_*(H)}$ is defined as 
\begin{gather*}
    CF_*(H)=\Big\{\sum_{z\in \widetilde{\mathcal{P}}(H)}a_z\cdot z \ | \ a_z\in \mathbb{Q},\forall C\in \mathbb{R}, \sharp \{z\in \widetilde{\mathcal{P}}(H) \ | \ a_z\neq 0, A_H(z)>0\}<\infty \Big\}.
\end{gather*}
The boundary operator ${d_F:CF_*(H)\rightarrow CF_{*-1}(H)}$ is defined as follows:
\begin{gather*}
    d_F(z)=\sum_{w\in \widetilde{\mathcal{P}}(H)}n(z,w)\cdot w
\end{gather*}
This coefficient ${n(z,w)\in \mathbb{Q}}$ is the number of solutions of the following Floer equation modulo the natural $\mathbb{R}$-action (see \cite{FO,LT}). We choose an almost complex structure $J$ on the tangent bundle ${TM}$.
\begin{gather*}
    z=[v_-,x_-], \ w=[v_+,x_+]  \\
    u:\mathbb{R}\times S^1\longrightarrow M \\
    \partial _su(s,t)+J\big(\partial_tu(s,t)-X_{H_t}(u(s,t))\big)=0  \\
    \lim_{s\rightarrow \pm\infty}u(s,t)=x_{\pm}(t)  \\
    [v_-\sharp u,x_+]=[v_+,x_+]
\end{gather*}
The Floer homology ${HF_*(H)}$ is the homology of the chain complex ${(CF_*(H),d_F)}$. It is known that ${HF_*(H)}$ is isomorphic to the quantum homology of ${(M,\omega)}$. Let ${\Gamma_{(M,\omega)}}$ be an abelian group defined as follows:
\begin{gather*}
    \Gamma_{(M,\omega)}=\frac{\pi_2(M)}{\mathrm{ker}\omega \cap \mathrm{ker}c_1}
\end{gather*}
Here, ${\omega:\pi_2(M)\rightarrow \mathbb{R}}$ is the integration of the symplectic form $\omega$ and ${c_1:\pi_2(M)\rightarrow \mathbb{Z}}$ is the integration of the first Chern class $c_1$. The degree of ${u\in \Gamma_{(M,\omega)}}$ is ${-2c_1(u)}$. The Novikov ring ${\Lambda_{(M,\omega)}}$ is defined as the set of possibly infinite sum of ${\Gamma_{(M,\omega)}}$ with suitable convergence, i.e.,
\begin{gather*}
    \Lambda_{(M,\omega)}=\Big\{\sum_{u\in \Gamma_{(M,\omega)}}a_u\cdot u \ | \ a_u\in \mathbb{Q}, \forall C\in \mathbb{R}, \sharp \{u\in\Gamma_{(M,\omega)}, \ a_u\neq 0, \omega(u)<C\}<\infty \Big\}.
\end{gather*}
The quantum homology of ${(M,\omega)}$ is the singular homology with ${\Lambda_{(M,\omega)}}$ coeffient, i.e.,
\begin{gather*}
    QH_*(M,\omega)=H_*(M:\mathbb{Q})\otimes \Lambda_{(M,\omega)}.
\end{gather*}
It is known that there is a natural isomorphism (PSS-isomorphism) between ${QH_*(M,\omega)}$ and ${HF_*(H)}$ (\cite{PSS,FO,LT}).
\begin{gather*}
    \mathrm{PSS}:QH_*(M,\omega)\longrightarrow HF_{*-n}(H)
\end{gather*}
The PSS-isomorphism is used to define the spectral invariant of the Floer homology. For any nonzero chain ${c=\sum_{z\in \widetilde{\mathcal{P}}(H)} a_z\cdot z\in CF_*(H)}$, the set ${\{A_H(z) | a_z\neq 0\}}$ is bounded above and discrete. We denote its supremum by ${l(c)}$. If $c=0$, we define ${l(0)=-\infty}$.
\begin{gather*}
    l:CF_*(H)\longrightarrow \mathbb{R}\cup \{-\infty\} \\
    c\mapsto \sup \{A_H(z) | a_z\neq 0\}
\end{gather*}
Then, ${CF_*^{<a}(H)=\{c\in CF_*(H)|l(c)<a\}}$ is a subcomplex of ${(CF_*(H),d_F)}$ for any ${a\in \mathbb{R}}$. We denote the homology of ${(CF_*^{<a}(H),d_F)}$ by ${HF_*^{<a}(H)}$. The inclusion ${CF_*^{<a}(H)\rightarrow CF_*(H)}$ induces a natural map 
\begin{gather*}
    \iota_a:HF_*^{<a}(H)\longrightarrow HF_*(H).
\end{gather*}
For any ${\alpha\in QH_*(M,\omega)}$, the spectral invariant ${c(\alpha, H)}$ is defined as follows (\cite{Oh,Sc,V}):
\begin{gather*}
    c(\alpha,H)=\inf\{a\in \mathbb{R} \ | \ \mathrm{PSS}(\alpha)\in \mathrm{Im}(\iota_a)\}
\end{gather*}
Usher proved that ${c(\alpha,H)}$ is achieved by some cycle in ${\mathrm{PSS}(\alpha)\in HF_*(H)}$.
\begin{Lem}[\cite{U2}]
    Let $H\in C^{\infty}(S^1\times M)$ be a non-degenerate Hamiltonian function. For any ${\alpha \in QH_*(M,\omega)\backslash\{0\}}$, there is a cycle ${c\in CF_*(H)}$ such that
    \begin{gather*}
        l(c)=c(\alpha,H)  \\
        [c]=\mathrm{PSS}(\alpha)
    \end{gather*}
    holds.
\end{Lem}

Note that ${c(\alpha,H)}$ is defined for possibly degenerate Hamiltonian function $H$ because ${c(\alpha,H)}$ is continuous with respect to the Hofer norm ($=$$L^{\infty}$-norm) of the function $H$. 
\begin{gather*}
    |c(\alpha,H)–c(\alpha,K)|\le \int_0^1\max_{x\in M}|H(t,x)-K(t,x)|dt
\end{gather*}
This continuity enable us to extend ${c(\alpha,\cdot)}$ to continuous functions.
The spectral norm of $H$ is defined by the spectral invariant of $H$ and $\overline{H}$ with respect to the funcamental class ${[M]\in QH_{2n}(M,\omega)}$.
\begin{gather*}
    \rho(H)=c([M],H)+c([M],\overline{H})
\end{gather*}
It is known that ${\rho(H)\ge 0}$ holds and moreover, ${\rho(H)>0}$ holds if ${\phi_H\neq \mathrm{Id}}$. This is a consequence of the proof of the energy-capacity inequality between the Hofer-Zehnder capacity and the displacement energy (see \cite{U,S}). More precisely, we have the following estimate of ${\rho(H)}$. Let ${(\mathbb{R}^{2n},\omega_0)}$ be the standard symplectic space.
\begin{gather*}
    \omega_0(x_1,\cdots,x_n,y_1,\cdots,y_n)=\sum_i x_i\wedge y_i
\end{gather*}
We say that ${\phi\in \mathrm{Ham}(M,\omega)}$ displaces a symplectically embedded $r$-ball (${r>0}$) if there is a symplectic embedding $\iota$
\begin{gather*}
    \iota:\big(B(r),\omega_0|_{B(r)}\big)\longrightarrow (M,\omega) \\
    B(r)=\big\{(x_1,\cdots,x_n,y_1,\cdots,y_n)\in \mathbb{R}^{2n} \ \big| \ \sum_ix_i^2+\sum_iy_i^2 \le r^2\big\}
\end{gather*}
such that ${\phi(\iota(B(r)))\cap \iota(B(r))=\emptyset}$ holds. If $\phi_H$ displaces a symplectically embedded $r$-ball,
\begin{gather*}
    \rho(H)\ge \pi r^2
\end{gather*}
holds. The Poincare duality of the spectral invariant (see \cite{LZ}) implies that if $H$ is non-degenerate, 
\begin{gather*}
    c([M],\overline{H})=-\inf\{c(\beta,H) \ | \ \beta\in QH_0(M,\omega), \beta=[pt]+\sum_{b\in H_{2k}(M), k\ge 1}a_b\cdot b\}
\end{gather*}
holds. In particular, this equality and Lemma 1 implies that there exists ${z,w\in \widetilde{\mathcal{P}}(H)}$ such that 
\begin{gather*}
    \mu_{CZ}(z)=n  \\
    \mu_{CZ}(w)=-n  \\
    A_H(z)=c([M],H), A_H(w)=-c([M],\overline{H}) 
\end{gather*}
holds. In particular, if $\phi_H$ displaces a symplectically embedded $r$-ball, there exists ${z,w\in \widetilde{\mathcal{P}}(H)}$ such that
\begin{gather*}
    0\le \mu(z)\le 2n  \\
    -2n\le \mu(w) \le 0 \\
    A_H(z)-A_H(w)\ge \pi r^2
\end{gather*}
holds. The same conclusion holds for possibly degenerate Hamiltonian function $H$.
\begin{Lem}
    Let $H\in C^{\infty}(S^1\times M)$ be a possibly degenerate Hamiltonian function such that $\phi_H$ displaces a symplectically embedded $r$-ball. Then, there are capped periodic orbits ${z,w\in \widetilde{\mathcal{P}}(H)}$ such that 
    \begin{gather*}
        0\le \mu(z)\le 2n \\
        -2n\le \mu(w)\le 0 \\
        A_H(z)-A_H(w)\ge \pi r^2
    \end{gather*}
    holds.
\end{Lem}
First, we choose a family of Hamiltonian functions ${\{H_m\}_{m\in \mathbb{N}}}$ such that every $H_m$ is non-degenerate and ${H_m}$ converges to $H$ in $C^{\infty}$-topology. We assume that every $\phi_{H_m}$ displaces a symplectically embedded $r$-ball. Then, we can choose ${z_m=[u_m,x_m]}$, ${w_m=[v_m,y_m]}$ in ${\widetilde{\mathcal{P}}(H_m)}$ so that
\begin{gather*}
    0\le \mu(z_m)\le 2n, \ -2n\le \mu (w_m)\le 0  \\
    A_{H_m}(z_m)=c([M],H_m), \ A_{H_m}(w_m)=c([M],\overline{H_m}) \\
    A_{H_m}(z_m)-A_{H_m}(w_m)\ge \pi r^2
\end{gather*}
Without loss of generality, we assume that $x_m$ and $y_m$ converges to contractible periodic orbits of $H$.
\begin{gather*}
    x_m\rightarrow x, \ y_m\rightarrow y
\end{gather*}
Note that the continuity of the spectral invariant implies that ${A_{H_m}(z_m)}$ converges to ${c([M],H)}$ and ${A_{H_m}(w_m)}$ converges to ${c([M],\overline{H})}$. Let $C_m$ and $D_m$ be small cylinders ${[0,1]\times S^1\rightarrow M}$ which connect $x_m$ and $x$, $y_m$ and $y$ respectively.
\begin{gather*}
    C_m(0,t)=x_m(t), \ C_m(1,t)=x(t) \\
    D_m(0,t)=y_m(t), \ D_m(1,t)=y(t)
\end{gather*}
Here, ``small" means that $\{C_m(s,t)\}_{s\in [0,1]}$ and ${\{D_m(s,t)\}_{s\in [0,1]}}$ are shortest geodesics with respect to some Riemanninan metric. We define ${\widetilde{z}_m}$ and ${\widetilde{w}_m}$ in ${\widetilde{\mathcal{P}}(H)}$ as follows:
\begin{gather*}
    \widetilde{z}_m=[u_m\sharp C_m,x], \ \widetilde{w}_m=[v_m\sharp D_m,y]
\end{gather*}
Then, $\{A_{H}(\widetilde{z}_m)\}$ converges to ${c([M],H)}$ and ${A_{H}(\widetilde{w}_m)}$ converges to ${c([M],\overline{H})}$. Our monotonicity assumption implies that the differences are contained in the discrete subgroup of $\mathbb{R}$.
\begin{gather*}
    A_H(\widetilde{z}_m)-A_H(\widetilde{z}_{m'})\in \frac{1}{\kappa}\mathbb{Z}  \\
    A_H(\widetilde{w}_m)-A_H(\widetilde{w}_{m'})\in \frac{1}{\kappa}\mathbb{Z}
\end{gather*}
This implies that ${A_H(\widetilde{z}_m)}$ and ${A_H(\widetilde{w}_m)}$ stabilize for sufficiently large $m$. In other words, we can choose $N\in \mathbb{N}$ so that
\begin{gather*}
    A_H(\widetilde{z}_m)=A_H(\widetilde{z}_{m'}), \ A_H(\widetilde{w}_m)=A_H(\widetilde{w}_{m'})    
\end{gather*}
holds for any $m,m'\ge N$. Our monotonicity assumption implies that the mean index also stabilizes. So, 
\begin{gather*}
    \mu(\widetilde{z}_m)=\mu(\widetilde{z}_{m'}), \ \mu(\widetilde{w}_m)=\mu(\widetilde{w}_{m'})    
\end{gather*}
holds for any $m,m'\ge N$. This also implies that $\widetilde{z}_m=\widetilde{z}_{m'}$ and ${\widetilde{w}_m=\widetilde{w}_{m'}}$ hold. We denote these capped periodic orbits by $z$ and $w$.
\begin{gather*}
    z=\widetilde{z}_m, \ w=\widetilde{w}_m \ \ \ m\ge N
\end{gather*}
This $z$ and $w$ satisfies the following relations:
\begin{gather*}
    0\le \mu(z) \le 2n, \ -2n\le \mu(w)\le 0  \\
    A_H(z)=c([M],H), \ A_H(w)=c([M],\overline{H})  \\
    A_H(z)-A_H(w)\ge \pi r^2
\end{gather*}
So we proved the lemma.

\section{Proof of the main theorem}
In this section, we prove the main theorem. First, we prove the following proposition.
\begin{Prop}
Let $(M,\omega)$ be a closed symplectic manifold. We fix a sufficiently small $C^1$-open neighborhood of the identity ${\mathcal{U}\subset \mathrm{Ham}(M,\omega)}$. For any ${\psi \in \mathcal{U}}$, we can construct a $C^1$-small Hamiltonian isotopy between the identity and ${\psi}$.
\end{Prop}
\begin{Rem}
It is almost trivial that we can construct a $C^1$-small symplectic isotopy between the identity and $\psi$. The important point is that we can construct a Hamiltonian isotopy.
\end{Rem}
We apply the following correspondence between the set of symplectomorphisms which are $C^1$-close to the identity and the set of small closed $1$-forms on $M$ (see \cite{MS}). Let ${\sigma \in \Omega(M\times M)}$ be a symplectic form on ${M\times M}$.
\begin{gather*}
    \sigma=-\pi_1^*\omega+\pi_2^*\omega
\end{gather*}
Here, ${\pi_i:M\times M\rightarrow M}$ is the projection from the $i$-th factor
\begin{gather*}
    \pi_i(x,y)=\begin{cases}
        x & i=1  \\
        y & i=2 .
    \end{cases}
\end{gather*}
We denote the canonical Liouville $1$-form on the cotangent bundle ${T^*M}$ by $\lambda$.
\begin{gather*}
    \lambda(X)=w(d\pi_*(X)) \ \ \ X\in T_wT^*M
\end{gather*}
This $d\pi:TT^*M\rightarrow TM$ is the differentail of the natural projection ${\pi:T^*M\rightarrow M}$.  Then, ${d\lambda \in \Omega(T^*M)}$ is the canonical symplectic form on the cotangent bundle. Let ${M_0\subset T^*M}$ be the image of the zero section. We can construct a symplectomorphism $\Psi$ between an open neighborhood of the zero section ${\mathcal{N}(M_0)\subset T^*M}$ and an open neighborhood of the diagonal ${\mathcal{N}(\Delta)\subset M\times M}$. This symplectomorphism $\Psi$ satisfies the following properties.
\begin{gather*}
    \Psi:\mathcal{N}(\Delta)\longrightarrow \mathcal{N}(M_0) \\
    \Psi(\Delta)=M_0  \\
    \pi(\Psi(q,q))=q \in M_0
\end{gather*}
We fix a Hamiltonian diffeomorphism ${\psi\in \mathrm{Ham}(M,\omega)}$ which is sufficiently $C^1$-close to the identity. We define a closed $1$-form ${\sigma_{\psi}}$ by 
\begin{gather*}
    \sigma_{\psi}=\Psi(\mathrm{graph}(\psi)).
\end{gather*}
Next we fix a $1$-periodic Hamiltonian function ${H\in C^{\infty}(S^1\times M)}$ so that ${\phi_H=\psi}$ holds. We define a path of symplectomorphisms ${\{\psi^t\}_{0\le t\le 1}}$ as follows:
\begin{gather*}
    \{(q,\psi^t(q))\}_{q\in M}=\Psi^{-1}(t\sigma_{\psi})
\end{gather*}
In particular, ${\sigma_{\psi^t}=t\sigma_{\psi}}$ holds. Our purpose is to prove that ${\{\psi^t\}}$ is a Hamiltonian isotopy. It sufficies to prove that $\sigma_{\psi}$ is an exact $1$-form because the flux of the path ${\{\psi^t\}_{0\le t\le T}}$ is equal to ${T[\sigma_{\psi}]}$ for any ${T\in [0,1]}$. Let ${\gamma_t}$ be a loop of symplectomorphisms defined as follows:
\begin{gather*}
    \gamma_t=\begin{cases}
        \psi_{2t} & 0\le t\le \frac{1}{2} \\
        \phi_H^{2(1-t)} & \frac{1}{2}\le t \le 1
    \end{cases}
\end{gather*}
The flux of the loop ${\{\gamma^t\}}$ is equal to $[\sigma_{\psi}]\subset \Gamma$. Note that $\Gamma$ is a discrete subgroup of ${H^1(M:\mathbb{R})}$. This implies that ${[\sigma_{\phi}]=0}$ holds if ${\sigma_{\psi}}$ is a sufficiently small closed $1$-form on $M$. In particular, ${\{\psi^t\}}$ is a Hamiltonian path. So, we proved the proposition.

\begin{Rem}
This argument can be used to prove that $C^1$-topology is stronger than the topology induced from Hofer's metric (\cite{S3}).
\end{Rem}

Henceforth, we assume that ${(M,\omega)}$ is a $2n$-dimensional closed negatively monotone symplectic manifold. Let ${H\in C^{\infty}(S^1\times M)}$ be a non-degenerate Hamiltonian function such that ${\phi_H}$ displaces a symplectically embedded $r$-ball. By lemma 2, we can choose two capped periodic orbits ${\Bar{x},\Bar{y}\in \widetilde{\mathcal{P}}(H)}$ so that 
\begin{gather*}
   -n\le  \mu(\Bar{x})\le 0\\
   0\le  \mu(\Bar{y})\le n  \\
    A_H(\Bar{y})-A_H(\Bar{x})\ge \pi r^2
\end{gather*}
holds. The $k$-th power ${\Bar{x}^{(k)},\Bar{y}^{(k)}\in \widetilde{\mathcal{P}}(H^{(k)})}$ satisfies
\begin{gather*}
    \mu(\Bar{x}^{(k)})\le 0\le \mu(\Bar{y}^{(k)}) \\
    A_{H^{(k)}}(\Bar{y}^{(k)})-A_{H^{(k)}}(\Bar{x}^{(k)})=k\{A_H(\Bar{y})-A_H(\Bar{x})\}\ge k\pi r^2 .
    \end{gather*}
We choose ${A_k,B_k\in \pi_2(M)}$ so that
\begin{gather*}
    0\le \mu(\Bar{x}^{(k)}\sharp A_k)\le 2N \\
    0\le \mu(\Bar{y}^{(k)}\sharp B_k)\le 2N
\end{gather*}
holds. This ${N\in \mathbb{N}\cup \{+\infty\}}$ is the minimal Chern number of ${(M,\omega)}$. In other words, $N$ is the positive generator of the image of ${c_1:\pi_2(M)\rightarrow \mathbb{Z}}$. We assume that $B_k$ is trivial if 
\begin{gather*}
    0\le \mu(\Bar{y}^{(k)})\le 2N
\end{gather*}
holds. This implies that 
\begin{gather*}
    A_{H^{(k)}}(\Bar{y}^{(k)}\sharp B_k)-A_{H^{(k)}}(\Bar{x}^{(k)}\sharp A_k) \ge A_{H^{(k)}}(\Bar{y}^{(k)})-A_{H^{(k)}}(\Bar{x}^{(k)})  \\
    =k(A_H(\Bar{y})-A_H(\Bar{x}))\ge k\pi r^2
\end{gather*}
holds because 
\begin{gather*}
    c_1(A_k)\le 0, \  \omega(A_k)\ge 0
\end{gather*}
holds. Next, we assume that ${(\phi_H)^k\subset \mathcal{U}}$ holds where ${\mathcal{U}\subset \mathrm{Ham}(M,\omega)}$ is a $C^1$-small open neighborhood of the identity. Proposition 1 implies that we can connect ${(\phi_H)^k}$ and ${\mathrm{Id}}$ by a $C^1$-small Hamiltonian isotopy generated by a Hamiltonian function ${K\in C^{\infty}(S^1\times M)}$. In particular, 
\begin{gather*}
    \phi_{H^{(k)}\sharp K}=\mathrm{Id}
\end{gather*}
holds. By concatenating this small isotopy, ${\Bar{x}^{(k)}\sharp A_k}$ and ${\Bar{y}^{(k)}\sharp B_k}$ become capped periodic orbits ${z_k,w_k\in \widetilde{\mathcal{P}}(H^{(k)}\sharp K)}$ such that
\begin{gather*}
    -\epsilon\le \mu(z_k) \le 2N+\epsilon \\
    -\epsilon\le \mu(w_k) \le 2N+\epsilon \\
    A_{H^{(k)}\sharp K}(w_k)-A_{H^{(k)}\sharp K}(z_k)\ge k\pi r^2-\epsilon
\end{gather*}
holds for some ${\epsilon >0}$. ${\epsilon >0}$ is determined by ${\mathcal{U}}$. Note that ${\phi_{H^{(k)}\sharp K}=\mathrm{Id}}$ holds. We compare ${\widetilde{\mathcal{P}}(H^{(k)}\sharp K)}$ and ${\widetilde{\mathcal{P}}(0)=\Gamma_{(M,\omega)}}$. More generally, we can construct a one to one correspondence between ${\widetilde{\mathcal{P}}(G_1)}$ and ${\widetilde{\mathcal{P}}(G_2)}$ if ${\phi_{G_1}=\phi_{G_2}}$ holds (see \cite{S2}, section 3.1). We fix a Hamiltonian function $L$ so that
\begin{gather*}
    \phi_{G_2}^t=\phi_L^t\circ \phi_{G_1}^t
\end{gather*}
holds. Note that $L$ generates a Hamiltonian loop. Let ${\mathcal{L}(M)}$ be the space of contractible loops in $M$.
\begin{gather*}
    \mathcal{L}(M)=\{x:S^1\rightarrow M \ | \ x:\mathrm{contractible}\}
\end{gather*}
We define a covering ${\pi:\widetilde{\mathcal{L}}(M)\rightarrow \mathcal{L}(M)}$ as follows:
\begin{gather*}
    \widetilde{\mathcal{L}}(M)=\{(u,x) \ | \ x\in \mathcal{L}(M),u:D^2\rightarrow M, \partial u=x\}/\sim  \\
    (u,x)\sim (v,y) \Longleftrightarrow \begin{cases}
        x=y \\ \int_{S^1}(u\sharp \overline{v})^*\omega=0  \\ \int_{S^1}(u\sharp \overline{v})^*c_1=0
        \end{cases}
\end{gather*}
The Hamiltonian loop ${\{\phi_L^t\}}$ acts on ${\mathcal{L}(M)}$ as follows:
\begin{gather*}
    f:\mathcal{L}(M)\longrightarrow \mathcal{L}(M) \\
    f(x)(t)=\phi_L^t(x(t))
\end{gather*}
Note that ${f(\mathcal{P})(G_1)=\mathcal{P}(G_2)}$ holds. Let ${\widetilde{f}}$ be a covering transformation such that the following diagram is commutative.
$$
\begin{CD}
\widetilde{\mathcal{L}}(M)  @>\widetilde{f}>> \widetilde{\mathcal{L}}(M)  \\
@V\pi VV  @V\pi VV  \\
\mathcal{L}(M) @>f>> \mathcal{L}(M)
\end{CD}
$$
We see that ${\widetilde{f}}$ restricted to ${\widetilde{\mathcal{P}}(G_1)}$ makes shifts of the action functional and the mean index.

\begin{Lem}
    For any ${z,w\in \widetilde{\mathcal{P}}(G_1)}$, the following equations holds:
    \begin{gather*}
        A_{G_1}(z)-A_{G_1}(w)=A_{G_2}(\widetilde{f}(z))-A_{G_2}(\widetilde{f}(w)) \\
        \mu(z)-\mu(w)=\mu(\widetilde{f}(z))-\mu(\widetilde{f}(w))
    \end{gather*}
\end{Lem}
For the first equation, we prove that ${A_{G_1}(z)-A_{G_2}(\widetilde{f}(z))}$ does not depend on ${z\in \widetilde{\mathcal{L}}(M)}$. It suffices to prove that the differential
\begin{gather*}
    T_z\widetilde{\mathcal{L}}(M)\longrightarrow \mathbb{R}  \\
    X\mapsto D(A_{G_1}(z)-A_{G_2}(\widetilde{f}(z)))X
\end{gather*}
vanishes for any ${z=[u,x]}$. Note that ${T_z\widetilde{\mathcal{L}}(M)=\Gamma(x^*TM)}$ holds. So $X$ is a section of the vector bundle ${x^*TM\rightarrow S^1}$.
\begin{gather*}
    D(A_{G_1}(z))X=-\int_{S^1}\omega(X(t),\dot{x}(t))dt+\int_{S^1}d(G_1)_t\cdot X(t)dt
\end{gather*}
\begin{gather*}
    D(A_{G_2}(\widetilde{f}(z)))X  \\
    =-\int_{S^1}\omega(d\phi_L^t(X(t)),\frac{\partial}{\partial t}(\phi_L^t(x(t))))dt+\int_{S^1}d(G_2)_t \cdot (d\phi_L^t\cdot X(t))dt \\
    =-\int_{S^1}\omega(d\phi_L^t(X(t)), X_{L_t}+d\phi_L^t(\dot{x}(t)))dt \\
    +\int_{S^1}(dL_t+d(G_1)_t\circ d(\phi_L^t)^{-1})\circ (d\phi_L^t\cdot X(t))dt \\
    =-\int_{S^1}\Big\{\omega(X(t),\dot{x}(t))+dL_t\circ d\phi_L^t(X(t))\Big\}dt \\
    +\int_{S^1}\Big\{dL_t\circ d\phi_L^t(X(t))+d(G_1)_t\cdot X(t)\Big\}dt  \\
    =-\int_{S^1}\omega(X(t),\dot{x}(t))dt+\int_{S^1}d(G_1)_t\cdot X(t)dt \\
    =D(A_{G_1}(z))X
\end{gather*}
So, ${D(A_{G_1}(z)-A_{G_2}(\widetilde{f}(z)))}$ is zero and the difference ${A_{G_1}(z)-A_{G_2}(\widetilde{f}(z))}$ is a constant function on ${\widetilde{\mathcal{L}}(M)}$. Next, we prove that ${\mu(z)-\mu(\widetilde{f}(z))}$ does not depend on ${z\in \widetilde{\mathcal{P}}(G_1)}$. We fix ${z=[u,x]}$ and ${w=[v,y]}$ in ${\widetilde{\mathcal{P}}(G_1)}$. We choose a cylinder ${C:[0,1]\times S^1\rightarrow M}$ which connects $x$ and $y$ and ${w=[u\sharp C,y]}$ holds. Let ${D:[0,1]\times S^1\rightarrow M}$ be a cylinder defined as follows:
\begin{gather*}
    D(s,t)=\phi_L^t(C(s,t))
\end{gather*}
We also fix a trivialization of the symplectic vector bundle ${C^*TM}$. Then $x$ and $y$ determines two paths of symplectomorphisms on the symplectic vecor space ${(\mathbb{R}^{2n},\omega_0)}$. ${\mu(z)-\mu(w)}$ is a difference of the mean index of these two paths. Similary, ${\mu(\widetilde{f}(z))-\mu(\widetilde{f}(w))}$ is caluculated by fixing a trivialization of ${D^*TM}$. One such trivialization is obtained by the trivialization of ${C^*TM}$ and the Hamiltonian loop ${\{\phi_L^t\}}$. So the equality
\begin{gather*}
    \mu(z)-\mu(w)=\mu(\widetilde{f}(z))-\mu(\widetilde{f}(w))
\end{gather*}
holds and we proved Lemma 3.

Lemma 3 and ${\phi_{H^{(k)}\sharp K}=\mathrm{Id}}$ implies that there is a transformation
\begin{gather*}
    \widetilde{f}:\widetilde{\mathcal{P}}(H^{(k)}\sharp K)\longrightarrow \widetilde{\mathcal{P}}(0)=\Gamma_{(M,\omega)}
\end{gather*}
such that
\begin{gather*}
    A_{H^{(k)}\sharp K}(z)-A_{H^{(k)}\sharp K}(w)=A_0(\widetilde{f}(z))-A_0(\widetilde{f}(w))  \\
   =\kappa(\mu(\widetilde{f}(z)-\mu(\widetilde{f}(w)))=\kappa (\mu(z)-\mu(w))
\end{gather*}
holds for any capped periodic orbits ${z,w\in \widetilde{\mathcal{P}}(H^{(k)}\sharp K)}$. However, this is a contradiction for sufficiently large $k\in \mathbb{N}$ because
\begin{gather*}
    |\mu(z_k)-\mu(w_k)|\le 2N+2\epsilon  \\
    |A_{H^{(k)}\sharp K}(z_k)-A_{H^{(k)}\sharp K}(w_k)|\ge k\pi r^2-\epsilon\rightarrow +\infty
\end{gather*}
holds. So, we proved Theorem 1 (1). 
Next, we prove Theorem 1 (2). We need the following Proposition. 
\begin{Prop}[\cite{A} Proposition 1.1]
    Let ${(M,\omega)}$ be a negatively monotone symplectic manifold whose Euler number is not zero or ${\pi_1(M)}$ has finite center. Then, the flux group $\Gamma$ is trivial.
\end{Prop}
We fix a symplectomorphism ${\psi\in \mathrm{Symp}_0(M,\omega)\backslash\{\mathrm{Id}\}}$. We choose a path of symplectomorphisms ${\{\psi^t\}}$ which connects ${\mathrm{Id}}$ and ${\psi}$ (${\psi^0=\mathrm{Id}}$, ${\psi^1=\psi}$). We extend ${\{\psi^t\}_{0\le t\le 1}}$ to an isotopy ${\{\psi^t\}_{t\in \mathbb{R}}}$ periodically. Let ${\mathcal{V}\subset \mathrm{Symp}_0(M,\omega)}$ be a sufficiently $C^1$-small open neighborhood of the identity so that any element in $\mathcal{V}$ can be connected to ${\mathrm{Id}}$ by a $C^1$-small isotopy of symplectomorphisms. Assume that ${\psi^k\in \mathcal{V}}$ holds. We connect $\mathrm{Id}$ and ${\psi^k}$ by a $C^1$-small isotopy of symplectomorphisms ${\{\phi^t\}}$. Let $\{\gamma^t\}$ be a loop of symplectomorphisms defined as follows:
\begin{gather*}
    \gamma^t=\begin{cases}
        \psi^{2kt} & 0\le t\le \frac{1}{2} \\
        \phi^{2(1-t)} & \frac{1}{2}\le t \le 1
    \end{cases}
\end{gather*}
Proposition 2 implies that 
\begin{gather*}
   0= [\mathrm{Flux}(\{\gamma^t\})]=k[\mathrm{Flux}(\{\psi^t\}_{0\le t\le 1})]-[\mathrm{Flux}(\{\phi^t\})]
\end{gather*}
holds. Note that ${[\mathrm{Flux}(\{\phi^t\})]\in H^1(M:\mathbb{R})}$ is contained in a small neighborhood of $0$ because ${\{\phi^t\}}$ is a $C^1$-small isotopy of symplectomorphisms. So ${[\mathrm{Flux}(\{\psi^t\}_{0\le t\le 1})]}$ is zero if $k$ is sufficiently large. In particular, ${\psi}$ is a Hamiltonian diffeomorphism. This implies that Theorem 1 (2) follows from Theorem 1 (1).

\section*{Acknowledgements}
  The author is supported by the Research Fellowships of Japan Society for the Promotion of Science for Young Scientists. He gratefully acknowledges his teacher Kaoru Ono and Manabu Akaho for continuous supports.

\end{document}